\documentclass[12pt,reqno]{amsart}
\usepackage{amsmath,amssymb}
\author{Theodore Voronov}
\address{Moscow State University}
\email{theodore@mech.math.msu.su}
\thanks{The research was supported in part by the Russian Foundation for
  Basic Research, under grant No 96-01-01440.}
\title[Variational complex]{Complex Generated by Variational Derivatives.
  Lagrangian Formalism of Infinite Order and a Generalized Stokes'
  Formula.}

\date{}
  \newtheorem{theorem}{Theorem}
  \newtheorem*{lemma}{Lemma}
\DeclareMathSymbol\leqslant{\mathrel}{AMSa}{"36} 
\DeclareMathSymbol\geqslant{\mathrel}{AMSa}{"3E} 
\DeclareMathOperator{\Supp}{Supp}

\renewcommand{\le}{\leqslant}
  \newcommand{\G}{\Gamma}
  \newcommand{\lra}{\longrightarrow}
  \newcommand{\field}{\mathbb}
  \newcommand{\At}{{\tilde A}}
  \newcommand{\Fot}{{\tilde F_1}}
  \newcommand{\Flt}{{\tilde F_l}}
  \newcommand{\dder}[3]{{\frac{\partial^{#1}}%
                        {\partial{#2}\dots\partial{#3}}}}
  \newcommand{\der}[2]{{\frac{\partial{#1}}%
                             {\partial{#2}}}}
  \renewcommand{\d}{{\,\,}{\bar{\smash{\!\!d}}}}
  \newcommand{\vder}[2]{{\frac{\delta{#1}}{\delta{#2}}}}
  \newcommand{\D}{\stackrel{\d}{\lra}}

\begin{document}
\maketitle
{\bf 1. The main construction.}
  Let $M^{n|m}$ be a smooth  (super)manifold.
  Consider the Lagrangians of $r|s$-dimensional paths, i.e. the maps
  $\G:I^{r|s}\lra M^{n|m}$, where $I^{r|s}:=I^r\times {\field R}^{0|s}$
  denotes the  $r|s$-cube.
  Denote the space of these Lagrangians by $\Phi^{r|s}$.
  The corresponding functionals of action are of the form
  $$
   S[\G]=\int_{I^{r|s}}\!Dt\,L(x(t),\dot x(t),\dots).
  $$
  Here and below
  we omit inessential common signs before the integrals.
  We assume that the Lagrangians may depend on derivatives
  of arbitrary (finite) order. If a Lagrangian depends on the derivatives
  of order $\le k$, then it is well known that the variation of the action
  in general includes the derivatives of order $\le 2k$.
  Here is the formula for the variation
  (cf.~\cite[p.~692]{sovgem}):
  $$
  \delta S=\int_{I^{r|s}}\!Dt\,\delta x^A(t)\,\vder{L}{x^A}
                                   (x(t),\dot x(t),\dots),
  $$
  where we introduced the notation
  $$
  \vder{L}{x^A}:=
  \sum_{l=0}^{\infty}(-1)^l(-1)^{\At(\Fot+\Flt)}\dder{l}{t^{F_1}}{t^{F_l}}
  \,\der{L}{x_{,F_1\dots F_l}^{\,A}}.
  $$
  (Mind that for the Lagrangians of arbitrary order  the ``total
  derivatives'' w.r.t. the time variables are well-defined.)
  For each Lagrangian $L$  we define a new Lagrangian $\d L$,
  which determines the action for the paths of dimension advanced by one:
  \begin{equation}\label{dif}
  \d L:= \dot x^A_{r+1} \vder{L}{x^A}.
  \end{equation}
  \begin{theorem}\label{chain}
    The definition $\d$ is coordinate independent. The identity
  $\d\,^2=0$ holds and the operation $\d$ makes the space
  of all Lagrangians a cochain complex:
  \begin{equation}\label{com}
  \dots\lra\Phi^{r|s}\D\Phi^{r+1|s}\D\Phi^{r+2|s}\D\dots
  \end{equation}
  \end{theorem}
  \begin{proof}
  The variational derivatives transform under change of coordinates as the
  components of a covector. This implies
  coordinate independence of $\d$.
  Consider the square of $\d$.
  A straightforward calculation of $\d\d L$ in the general case is
  unavailable,
  since it would require
  the explicit expansion of all multiple total derivatives
  in the formula~(\ref{dif}).
  For the Lagrangians of first order the formal expression
  for $\d\,^2$ contains 36 terms ! They do cancel, indeed,
  but checking this fact (the author has done it) is far from pleasant.
  Let us prove the identity $\d\,^2=0$ geometrically.
  It is convenient to introduce  ``weighted'' paths and consider functionals
  of the form
  $$
   S[\G,g]=\int_{I^{r|s}}\!Dt\,g(t)L(x(t),\dot x(t),\dots),
  $$
  where the function $g$ is compactly supported in the cube $I^{r|s}$.
  The following lemma holds.
  \begin{lemma}\label{var} For a variation of the map $\G$
    $$
     \delta S=0 \quad\text{for}\quad \forall g\, \forall\delta x(t):
     \Supp \delta x\,\cap\,\Supp dg=\varnothing \Longleftrightarrow
     \vder{L}{x^A}(x(t),\dot x(t),\dots)=0.
    $$
  \end{lemma}
  \begin{proof} For the variations $\delta x^A$ subjected to the above
  restriction it is not difficult to check
  that $\delta S=\int \!Dt\,g(t)\delta x^A\, \delta L/\delta x^A$.
  Suppose $\delta S=0$. Take $g$ supported inside a small ball,
  identically $1$ near its center, and take $\delta
  x^A(t)$ as a product of a delta-like function in the even variables
  by an arbitrary function of the odd variables. Passing to the limit
  and making use of the nondegeneracy of the Berezin integral,
  we obtain the equality
  $(\delta L/\delta x^A)(t_0,\tau)\equiv0$,
  where $t_0\in I^r$ is arbitrary (the ball's center).
  \end{proof}
  Consider the action $S^*$ that corresponds to the Lagrangian $\d L$.
  By definition of the variational derivative,
  $$
  \der{x^A}{t^{r+1}}\vder{L}{x^A}(x(t),\dots)=\der{}{t^{r+1}}(L)+
  \,\text{derivatives w.r.t.} \,\,t^F, F\neq r+1.
  $$
  Thus, the coefficient of $g^*(t^*)$ in the integral is a ``total
  divergence''. So,  $S^*$ can be rewritten as the integral of $\pm(\partial
  g^*/\partial t^{F^*}) h^{F^*}(t^*)$ (analog of the integral over boundary).
  Then for the variation  $\delta x^A(t^*)$
  satisfying the condition $\Supp \delta x\,\cap\,dg^*=\varnothing$, as it
  is easy to see, the equation $\delta S^*=0$ always holds, since the
  support of $\delta h^{F^*}$ is less or equal to the support of
  $\delta x^A$. Applying Lemma~\ref{var}, we obtain $\delta\,\d L/
  \delta x^A=0$ identically, and thus $\d\d L=0$.
  \end{proof}

  The cohomology of the cochain complex~(\ref{com}) is an invariant
  of a smooth manifold  $M$.

{\bf 2. Generalized Stokes' Formula.} Suppose we have a family of $D$-paths
  $\G_s$ (where $s\in[0,1]$) that is constant outside of a compact
  in $I^{D}$, $D=r|s$.
  We consider it as a path $\G^*$, of dimension $D+1$.
  As an immediate corollary of the above consideration we obtain the
  following
  \begin{theorem} For any Lagrangian
  \begin{equation}\label{stokes}
     S[\G_1]-S[\G_0]= \int_{\G^*} \d L.
  \end{equation}
  \end{theorem}
  If the Lagrangian is such that the action is independent of
  para\-metri\-za\-tion (with an appropriate restriction on orientation,
  see~\cite{git}), then it is possible to
  integrate it over singular manifolds with boundary, and
  a complete analog of the Stokes' Formula is valid.
  We see that it is no exclusive property of the
  differential forms. The formula~(\ref{stokes}) provides explanation of
  the identity $\d\,^2=0$ as a manifestation of the general fact
  that the  ``boundary of boundary is zero''.

{\bf 3. Application to supermanifold forms.}
  Definition~(\ref{dif}) for the first time was suggested by the author
  in the following context~\cite{first}: one had to define the right analog
  of the de Rham complex for a supermanifold that is not an even manifold.
  It was suggested to consider as forms  the ``covariant'' Lagrangians
  of the first order, satisfying an extra condition that provided
  the independence of variation of the accelerations (see~\cite{first,git}).
  For the even case it was proved
  that the de~Rham--Cartan theory was reproduced. In  general situation
  the theory is nontrivial.
  The indicated condition is equivalent to a system
  of PDE~\cite[p.~58]{git}, for which the odd-odd part gives the equations
  well known in the connection with Radon-like
  integral transforms~\cite{john, ggg}.
  Let us introduce a filtration of Lagrangians by the
  order of derivatives.
  Then our system of equations receive a ``homological''
  interpretation as the condition for the elements under consideration
  to keep their filtration under the action of the differential $\d$.
  The higher analogs of this
  system are obvious in such interpretation.
  The results of the present note provide substantial simplification
  of the proofs from~\cite{git}.
  Recently a ``second half'' of the complex of supermanifold forms
  was constructed, with the use of the ``dual'' language~\cite{dual}.
  I hope that the combination of the methods of~\cite{dual} with those
  of the
  current paper will lead to new advances.

{\bf 4. Application to the inverse problem of the calculus of variations.}
  According to Theorem~\ref{chain}, in order to check,
  if some set of functions $f_A$ coincide with variational derivatives
  of a Lagrangian of  $r|s$-path, it is necessary to calculate
  $\d (\dot x_{r+1}^A f_A)$ and check that it is zero.
  (Note that the cross-sections of bundles also can be
  treated as ``surfaces'' or ``paths''.) The inverse problem of the
  calculus of variations  has been studied in~\cite{vin,tak} and others.
  It would be interesting to compare the approach of these papers
  with our complex~(\ref{com}).

  As a conclusion I wish to notice that
  the constructions of differentials known in homology theory
  fall into two or three main types (combinatorial differentials
  and E.~Cartan or Koszul type differentials).
  The ``variational'' differentials, studied in~\cite{first,git,dual}
  and in the present note, are likely to supply yet another type.

\end{document}